\documentclass[12pt]{article} 
\usepackage{amsmath}
\usepackage{amsfonts,amssymb}
\usepackage{fullpage}

\newcommand{\cluster}{\mathcal{C}}
\newcommand{\clustern}{\mathcal{C}_n}
\newcommand{\ball}{\mathcal{B}}
\newcommand{\balln}{\mathcal{B}_n}

\def\squareforqed{\hbox{\rlap{$\sqcap$}$\sqcup$}}
\def\qed{\ifmmode\else\unskip\quad\fi\squareforqed}
\def\smartqed{\def\qed{\ifmmode\squareforqed\else{\unskip\nobreak\hfil
\penalty50\hskip1em\null\nobreak\hfil\squareforqed
\parfillskip=0pt\finalhyphendemerits=0\endgraf}\fi}}

\newtheorem{thm}{Theorem}
\newtheorem{lem}[thm]{Lemma}

\newcommand{\aston}{\xrightarrow[n\to\infty]{a.s.}}
\begin{document}
\title{Long-range percolation on the hierarchical lattice}
\author{Vyacheslav Koval$^{1}$, Ronald Meester$^{2}$ and Pieter Trapman$^{3}$}
\footnotetext[1]{Utrecht University, Department of Mathematics, P.O.~Box 80010, 3508 TA Utrecht, The Netherlands,}
\footnotetext[2]{VU University Amsterdam, Department of Mathematics, De Boelelaan 1081a, 1081 HV Amsterdam, The Netherlands,}
\footnotetext[3]{Stockholm University, Department of Mathematics, 106 91 Stockholm, Sweden.}
\date{\today}
\maketitle
\begin{abstract}
We study long-range percolation on the hierarchical lattice of order $N$, where any edge of length $k$ is present with probability $p_k=1-\exp(-\beta^{-k} \alpha)$,
independently of all other edges. For
fixed $\beta$, we show that the critical value $\alpha_c(\beta)$ is non-trivial if and only if $N < \beta < N^2$. Furthermore, we show uniqueness of the infinite component and continuity of the percolation probability and of $\alpha_c(\beta)$ as a function of $\beta$. This means that the phase diagram of this model is well understood.
\end{abstract}

\section{Introduction and main results}

The use of percolation theory in statistical physics has long been recognized. The study of long-range percolation on $\mathbb{Z}^d$ goes back to \cite{schulman83} and led to a series of interesting problems and results \cite{aizenman86,newmanChuck86,berger02,biskup04,trapman10}; see \cite[Section 2]{biskup09} for an extensive overview. In \cite{dawson07} asymptotic long-range percolation is studied on the hierarchical lattice $\Omega_N$ (to be defined below) for $N\to\infty$. The contact process on $\Omega_N$ for fixed $N$ has been studied in \cite{AtSw09}.
 
In this paper we study the case of finite $N$. Long range percolation on the hierarchical lattice is quite different from the usual lattice: classical methods break down and results are different. The results and methods in this paper could appeal to both mathematicians and physicists.    

For an integer $N\ge 2$, we define the set
\[
\Omega_N := \{\mathbf x=(x_1,x_2,\ldots):x_i\in\{0, 1, \ldots,N-1\},\sum_i x_i<\infty\},
\]
and define a metric on it by
\[
d(\mathbf x,\mathbf y) =  \begin{cases}
  0 &\text{if $\mathbf x=\mathbf y$,}\\
  \max\{i:x_i\neq y_i\} &\text{if $\mathbf x\neq\mathbf y$.}
 \end{cases}
\]
The pair $(\Omega_N, d)$ is called the \emph{hierarchical lattice of order $N$}.
\begin{center}
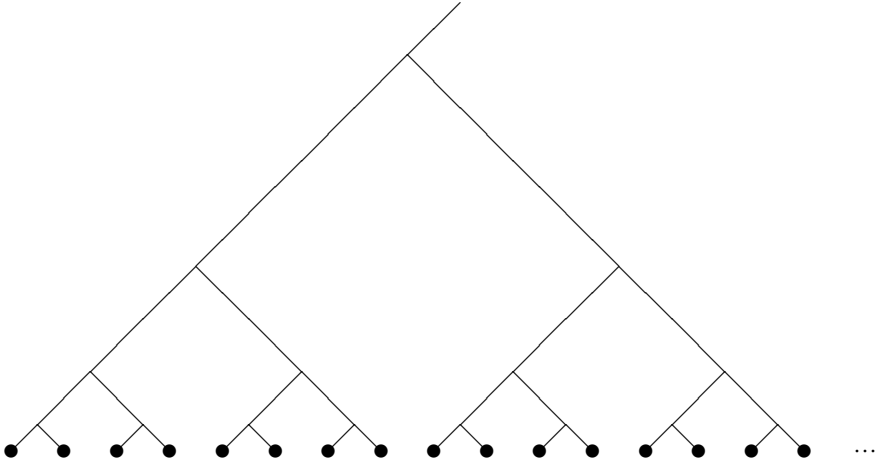
\begin{figure}
\begin{picture}(340,200)
\put(0,0){\circle*{5}}
\put(20,0){\circle*{5}}
\put(40,0){\circle*{5}}
\put(60,0){\circle*{5}}
\put(80,0){\circle*{5}}
\put(100,0){\circle*{5}}
\put(120,0){\circle*{5}}
\put(140,0){\circle*{5}}
\put(160,0){\circle*{5}}
\put(180,0){\circle*{5}}
\put(200,0){\circle*{5}}
\put(220,0){\circle*{5}}
\put(240,0){\circle*{5}}
\put(260,0){\circle*{5}}
\put(280,0){\circle*{5}}
\put(300,0){\circle*{5}}
\put(320,0){\circle*{1}}
\put(323,0){\circle*{1}}
\put(326,0){\circle*{1}}
\put(0,0){\line(1,1){170}}
\put(20,0){\line(-1,1){10}}
\put(40,0){\line(1,1){10}}
\put(60,0){\line(-1,1){30}}
\put(80,0){\line(1,1){30}}
\put(100,0){\line(-1,1){10}}
\put(120,0){\line(1,1){10}}
\put(140,0){\line(-1,1){70}}
\put(160,0){\line(1,1){70}}
\put(180,0){\line(-1,1){10}}
\put(200,0){\line(1,1){10}}
\put(220,0){\line(-1,1){30}}
\put(240,0){\line(1,1){30}}
\put(260,0){\line(-1,1){10}}
\put(280,0){\line(1,1){10}}
\put(300,0){\line(-1,1){150}}
\end{picture}
\caption{Hierarchical lattice (the ``leaves'')  of order $2$ with the metric generating tree attached.}
\label{fig:tree}
\end{figure}
\end{center}
One can think of the vertices in the hierarchical lattice as the leaves of a regular tree without a root, see Figure \ref{fig:tree}. The metric $d$ can then be interpreted as the number of generations (levels) till the ``most recent common ancestor'' of two vertices.
Let $\mathbb{N}$ be the non-negative integers, including $0$ and $\mathbb{N}_+ := \mathbb{N} \setminus \{0\}$.
The set $\Omega_N$ is countable, and we can introduce a natural labeling of its vertices via the map
$f: \Omega_N \to \mathbb{N}$ given by
\[
f:(x_1,x_2,\ldots) = \sum_{i=1}^{\infty} x_i N^{i-1}.
\]
We will sometimes abuse notation and write $n$ for $f^{-1}(n) \in \Omega_N$.

The metric space $(\Omega_N,d)$ satisfies the strengthened version of the triangle inequality
\[
d(x,y)\le \max (d(x,z),d(z,y))
\]
for any triple $x,y,z \in \Omega_N$. Such spaces are called ultrametric (sometimes non-Archimedean) \cite{vanrooij78}. For $x \in \Omega_N$, define $\mathcal B_r(x)$ to be the ball of radius $r$ around $x$.
Several important geometrical properties follow from the definition of the space $(\Omega_N, d)$ and its ultrametricity:
\begin{enumerate}
\item $\mathcal B_r(x)$ contains $N^r$ vertices for any $x$;
\item for every $x\in \Omega_N$ there are $(N-1)N^{k-1}$ vertices at distance $k$ from it;
\item if $y\in \mathcal B_r(x)$ then $\mathcal B_r(x)=\mathcal B_r(y)$;
\item as a consequence of the previous property, for all $x$, $y$ and $r$ we either have 
$\mathcal B_r(x)=\mathcal B_r(y)$ or $\mathcal B_r(x)\cap\mathcal B_r(y)=\emptyset$.
\end{enumerate}

Now consider long-range percolation on $\Omega_N$. Every pair of vertices $(x,y)\in\Omega_N \times \Omega_N$ is (independently of all other edges) connected by a single edge with probability 
$$
p_k = 1-\exp\left(-\frac{\alpha}{\beta^k}\right),
$$ 
where $k=d(x,y)$ and where $0 \leq \alpha <\infty$ and $0<\beta<\infty$ are the parameters of the model.
The vertices $x \in \Omega_N$ and $y \in \Omega_N$ are in the same connected component if there exists a \textit{path} from $x$ to $y$, that is, if there exists a finite sequence $x=x_0,x_1, \ldots, x_n=y$ of vertices such that every pair $(x_{i-1},x_{i})$ of points with $1 \leq i \leq n$ shares an edge. The edges are not directed.

We denote the size of a set $S$ of vertices by $|S|$. The connected component (also called ``cluster'') containing the vertex $x$ is denoted by $\cluster(x)$. Since $|\cluster(x)|$ has the same distribution for every $x\in\Omega_N$ we may study $|\cluster(0)|$ instead of $|\cluster(x)|$. 

Let $\mathbb{P}_{\alpha,\beta}$ be the probability measure governing this percolation process (on the appropriate probability space and 
sigma-algebra) and $\mathbb{E}_{\alpha,\beta}$ the corresponding expectation operator. When no confusion is possible, we 
omit the subscripts $\alpha$ and $\beta$. Denote 
\begin{equation}\label{thetadef}
\theta(\alpha,\beta):=\mathbb{P}_{\alpha,\beta}\left(|\cluster(0)|=\infty\right).
\end{equation}
It follows from a standard coupling argument 
that $\theta(\alpha,\beta)$ is non-decreasing in $\alpha$ for any given $\beta$. 
Therefore, it is reasonable to define
\[
\alpha_{c}(\beta):=\inf\{\alpha \geq 0:\theta(\alpha,\beta) > 0\}.
\]

Throughout the paper we use the following notation. For a set $S$ of vertices, let  $\overline{S} := \Omega_N\setminus S$ denote its complement. The set $\clustern(x)$ is the cluster of vertices that are connected to the origin by a path that uses only vertices inside $\balln(x)$. For disjoint sets $S_1,S_2 \subset \Omega_N$, the event that at least one edge connects a vertex in $S_1$ with a vertex in $S_2$ is denoted by $S_1\leftrightarrow S_2$. The notation $S_1\not\leftrightarrow S_2$ denotes the event that such an edge does not exist. Let $\cluster_n^m(x)$ be the largest cluster in $\balln(x)$; if more such clusters exist, $\cluster_n^m(x)$ is defined to be one of them, chosen uniformly among all possible candidates. In any case,
\begin{equation*}\label{defmaxclus}
|\cluster_n^m(x)| = \max_{y \in \balln(x)}|\clustern(y)|.
\end{equation*}

\begin{thm}\label{theorem1}((non)-triviality of the phase transition)
\begin{itemize}
\item[(a)] $\alpha_{c}(\beta)=0$ for $\beta\le N$,
\item[(b)] $0<\alpha_{c}(\beta)<\infty$ for $N<\beta<N^2$,
\item[(c)] $\alpha_{c}(\beta)=\infty$ for $\beta\ge N^2$.
\end{itemize}
\end{thm}
\begin{thm}
\label{theorem1a} (uniqueness of the infinite cluster)
There is a.s.\ at most one infinite cluster, for any value of $\alpha$ and $\beta$.
\end{thm}

\begin{thm}\label{theorem2}(continuity of $\theta$)
The percolation function $\theta(\alpha,\beta)$ is continuous whenever $\alpha>0$.  
\end{thm}

\begin{thm}\label{theorem4}
(continuity of $\alpha_{c}(\beta)$)
The critical value $\alpha_{c}(\beta)$ is continuous for $\beta \in (0,N^2)$ and strictly increasing for $\beta \in [N,N^2)$. Finally, $\alpha_{c}(\beta)\nearrow\infty$ for $\beta\nearrow N^2$.  
\end{thm}

In order to prove Theorems \ref{theorem2} and \ref{theorem4}, we need the following result, which is interesting in its own right.

\begin{thm}\label{theorem3} (size of large components)
If $\alpha$ and $\beta$ are such that $\theta := \theta(\alpha,\beta)>0$,
then for every $\varepsilon>0$,
\begin{equation}\label{fractionlemma}
\lim_{k \to \infty }\mathbb{P}{}_{\alpha,\beta}\left(|\mathcal{C}^m_k(0)|> (\theta-\varepsilon) N^k\right)=1.
\end{equation}
\end{thm}

In the next two sections we prove Theorem \ref{theorem1} and Theorem \ref{theorem1a} respectively. After that we prove the remaining results, and we end
with a discussion about possible generalisations.

\section{Proof of Theorem \ref{theorem1}}

\noindent
{\em Proof of (a).}
Denote by $E_k$ the event that the origin shares an edge with at least one vertex at distance $k$. Then
\[
\mathbb{P}(E_k) = 1-\exp\left(-\frac{\alpha}{\beta^k}(N-1)N^{k-1}\right)
\]
and the events $(E_k)_{k \geq 1}$ are independent. 
It is easy to see that if $\beta\le N$ then $\sum_{k=1}^\infty\mathbb{P}(E_k)$ diverges for any $\alpha>0$. Therefore, by the second Borel-Cantelli lemma, infinitely many of the events $E_k$ occur with probability $1$ and $\theta(\alpha,\beta)=1$, for any $\alpha>0$ and $0 < \beta \leq N$. This implies $\alpha_{c}(\beta)=0$ for $0 < \beta \leq N$.
\hfill{$\qed$}

\medskip\noindent
{\em Proof of (c).}
By monotonicity it suffices to prove that $\alpha_{c}(N^2)=\infty$, so we now take $\beta = N^2$.
A straightforward computation shows that for every $j$,
\begin{equation}
\mathbb{P}(\ball_j(v) \leftrightarrow \overline{\ball_j(v)}) = 1-\exp\left(-\alpha N^j \frac{(N-1)}{N^2}\sum_{k=1}^{\infty}\frac{N^{j+k-1}}{N^{2(j+k-1)}}\right) = 1-\exp\left(-\frac{\alpha}{N}\right),\nonumber
\end{equation}
which is strictly less than 1, and independent of $j$. 
Let $n_0=0$ and $n_{i+1}=\inf\{n\ge n_i:\ball_{n_i}(0)\not\leftrightarrow\overline{\ball_n(0)}\}$.
Since 
$$\{\cluster(0) = \infty\} \subset \cap_{i=0}^{\infty} \{\ball_{n_{i}}(0) \leftrightarrow\overline{\ball_{n_{i}}(0)} \},$$ it is enough to prove that there a.s.\ 
exists $i$ such that $\ball_{n_{i}}(0)\not\leftrightarrow\overline{\ball_{n_{i}}(0)}$.
Because the events  $\{\ball_{n_{i}}(0) \leftrightarrow\overline{\ball_{n_{i}}(0)} \}$ are independent and all have the same probability strictly less than 1, the result follows.   
\hfill{$\qed$}

\medskip\noindent
{\em Proof of (b).}
The strict positivity of $\alpha_{c}(\beta)$ follows from the fact that
\begin{eqnarray*}
\sum_{k=1}^{\infty} (N-1)N^{k-1} p_k &=& \sum_{k=1}^{\infty} (N-1)N^{k-1}\left( 1-\exp(-\frac{\alpha}{\beta^k})\right) \leq  \\
&\leq &  \frac{\alpha (N-1)}{N} \sum_{k=1}^{\infty}\left(\frac{N}{\beta}\right)^k = \alpha  (N-1)\frac{1}{\beta-N},
\end{eqnarray*}
which can be made strictly smaller than 1 by choosing $\alpha$ small enough. Hence the expected number of edges from a given vertex is strictly smaller than 1, and by coupling with a subcritical branching process, the almost sure finiteness of the percolation cluster follows.

The second inequality is much more involved. Choose an integer $K$ and a real number $\eta$ such that
\begin{equation}\label{crucial}
\sqrt{\beta}<\eta\le \left(N^K-1\right)^{1/K};
\end{equation}
this is possible since $\sqrt{\beta}<N$.
We say that a ball of radius $nK$ is {\em good} if its largest connected component has size at least $\eta^{nK}$. Denote by $s_n$ the probability that a ball of radius $nK$ is good, that is,
\[
s_n := \mathbb{P}\left(|\cluster_{nK}^m(0)| \ge \eta^{nK}\right).
\]
By convention, we set $s_0=1$.

We say that a ball of radius $nK$ is {\em very good} if it is good and in addition its
largest component shares an edge with the largest component of the first (the one with the smallest index) good sub-ball in the same ball of radius $(n+1)K$. Note that according to this definition, the first good sub-ball of diameter $nK$ in a ball of radius
$(n+1)K$ is automatically very good.

Since $(N^K-1)\ge\eta^K$, $\ball_{(n+1)K}(0)$ will certainly be good if
\begin{itemize}
\item[(a)] it contains $N^K-1$ good sub-balls of radius $nK$, and 
\item[(b)] all these good sub-balls are very good. 
\end{itemize}
We next estimate the probability of the events in (a) and (b). 

Clearly, the number of good sub-balls of radius $nK$ in a ball of radius $(n+1)K$ has a binomial distribution with
parameters $N^K$ and $s_n$. Furthermore, given the collection of good sub-balls, the probability that the first such good sub-ball is very good is equal to 1, and the probability for any of the other good sub-balls to be very good is at least
$$
1-\varepsilon_n:=1-\exp\left(-\frac{\alpha}{\beta^{K}}\left(\frac{\eta^2}{\beta}\right)^{nK}\right),
$$
since the distance between two vertices in a ball of radius $(n+1)K$ is at most $(n+1)K$ and the largest component of a good sub-ball contains at least $\eta^{nK}$ vertices. 
We conclude that the number of very good sub-balls is stochastically larger than a random variable having a binomial distribution with
parameters $N^K$ and $s_n(1-\varepsilon_n)$. It follows that
\begin{equation}
\label{fractals}
s_{n+1}\ge\mathbb{P}\Bigl(Bin\left(N^K,s_n(1-\varepsilon_n)\right)\ge N^K-1\Bigr),
\end{equation}
where $Bin(n,p)$ denotes a random variable with a binomial distribution with parameters $n$ and $p$. Notice that since
$1/x > e^{-x}$ for all $x>0$, we have 
$$
\varepsilon_n = \left(\exp\left(-\left(\frac{\eta^2}{\beta}\right)^{nK}\right)\right)^{\frac{\alpha}{\beta^{K}}}
\leq \left(\left(\frac{\beta}{\eta^2}\right)^{\alpha K  \beta^{-K}}\right)^n,
$$
and hence for any $\delta >0$ we can find $\alpha$ so large that $\varepsilon_n \leq \delta^n$ for all $n$.

The expression in (\ref{fractals}) is very close to the iteration formulae in \cite{deme90}, and it pays to recall the setup from that paper first.

Let, for some given $m$, the map 
$$
\pi(p) =\mathbb{P}(Bin(m,p) \geq m-1)
$$
be defined. Define
$u_0=1$ and, for $n\in \mathbb{N}$
$$
u_{n+1} = \pi(p u_n).
$$
Writing $G_p(\cdot)$ for $\pi(p\cdot)$ we then obtain
\begin{equation}
\label{oud}
u_{n+1}= G_p(u_n).
\end{equation}
In \cite{deme90}, this iteration arises in the study of fractal percolation and it is shown that the limit $u= \lim_{n \to \infty} u_n$ always exists and is positive if and only if $p$ is so
large that the equation $G_p(x)=x$ has a positive solution. This result is very similar (and can be proved in the same way)
as the classical non-extinction criterium for branching procesess. When $p$ approaches 1, then the largest solution
of $G_p(x)=x$ also approaches 1, at essentially the same rate, and therefore also the limit $u$ approaches 1.

Now observe that (\ref{fractals}) can be rewritten as
\begin{equation}
\label {jawel}
s_{n+1}\geq G_{1-\varepsilon_n}(s_n).
\end{equation}

This is very similar to (\ref{oud}), the only difference being that the subscript of the iteration function depends on $n$
now. However, we already showed above that for $\alpha$ large enough, $\varepsilon_n$ goes down exponentially fast
at any given rate. It is not hard to believe that this implies that $s_n$ converges
to 1 exponentially fast, and we make this precise now.  

We have
$$
\mathbb{P}(Bin(n,p) \ge n-1)\ge 1-\binom{n}{2}(1-p)^2,
$$
and writing $C=\binom{N^K}{2}$ we arrive at the inequality
\begin{equation}\label{ito}
s_{n+1}\ge 1 - C(1-s_n+s_n\varepsilon_n)^2\ge 1 - C(1-s_n+\varepsilon_n)^2.
\end{equation}
Writing $\xi_n=1-s_n$ this gives
\begin{equation}
\label{ennu}
\xi_{n+1} \leq C(\xi_n + \varepsilon_n)^2.
\end{equation}
Choose first $\gamma$ so small that $4C \leq \gamma^{-1}$ and then $\alpha$ so large that $\varepsilon_n \leq \gamma^n$
and $\xi_1 \leq \gamma^2$.
In an inductive fashion, if $\xi_n\le\gamma^{n+1}$ then
\begin{equation}\label{largeprob}
\xi_{n+1}\le C(\xi_n+\varepsilon_n)^2
\le 4C(\gamma^{n+1})^2\le \gamma^{2n+1}\le \gamma^{n+2},
\end{equation}
which implies that $\xi_n\le\gamma^{n+1}$ for all $n$. Hence, for $\alpha$ large enough, $s_n$ converges to 1 exponentially
fast.

The exponential convergence of $s_n$ to 1 is not quite enough for our purposes. Indeed, $s_n$ represents
the probability that a ball of radius $nK$ contains a component of size at least $\eta^{nK}$, but this component does not
necessarily contain the origin. Therefore, we have to make one extra step. Let
$$
t_n:=\mathbb{P}\bigl(|\cluster_{nK}(0)|\ge \eta^{nK}\bigr).
$$
We claim that
\begin{equation}
\label{eindelijk}
t_{n+1} \geq t_n \times \mathbb{P}(Bin(N^K-1, s_n(1-\varepsilon_n))\geq N^K-2).
\end{equation}
To see this, we argue as before. If $|\cluster_{nK}(0)|\ge \eta^{nK}$, then $\ball_{nK}(0)$ will be the first good sub-ball in the 
derivation above. If this component is connected to at least $N^K-2$ other large components in $\ball_{(n+1)K}(0)$ as above, then the component of the origin in $\ball_{(n+1)K}$ is large enough, that is, has size at least $\eta^{(n+1)K}$.  From this, (\ref{eindelijk}) follows. Since a simple coupling gives that
$$
\mathbb{P}(Bin(N^K-1, s_n(1-\varepsilon_n))\geq N^K-2) \geq \mathbb{P}(Bin(N^K, s_n(1-\varepsilon_n))\geq N^K-1),
$$
and since the derivation above actually gives that the right hand side of this inequality converges to 
1 exponentially fast, it follows
that 
$$
\lim_{n \to \infty} t_n >0,
$$
which is enough to prove the result.
\hfill{$\qed$}

\medskip\noindent
\textbf{Remarks} \textbf{(I)} By the proof of the strict positivity of $\alpha_{c}(\beta)$ for $\beta>N$, we also may deduce that $\alpha_{c}(\beta)$ is not differentiable in $N$. Indeed, $\alpha_{c}(\beta) = 0$ for $\beta \leq N$ and  
since $$\sum_{k=1}^{\infty} (N-1)N^{k-1} p_k \leq \alpha  (N-1)\frac{1}{\beta-N}$$ for $\beta>N$, we have $\alpha_{c}(\beta) \geq \frac{\beta-N}{N-1}$ for all $\beta>N$. This in turn implies that for all $\beta>N$ we have  
$$\frac{\alpha_{c}(\beta)-\alpha_{c}(N)}{\beta-N} \geq \frac{1}{N-1} >0.$$
\medskip\noindent  \textbf{(II)} Since we may choose $\gamma$ arbitrary small in equation (\ref{largeprob}), we in fact have that for every $\varepsilon>0$, we can choose $\alpha$ so large, that $t_n>1-\varepsilon$. This implies that for every $\beta \in (N,N^2)$, we can choose $\alpha$ so large such that $\theta(\alpha, \beta)>1-\varepsilon$.

\section{Proof of Theorem \ref{theorem1a}}

We will use Theorem 0 from \cite{gandolfi92}:

\begin{thm} (Gandolfi et al.\ \cite{gandolfi92})
Consider long range percolation on $\mathbb{Z}^d$ with the properties
\begin{enumerate}
\item the model is translation-invariant, and
\item the model satisfies the positive finite energy condition.
\end{enumerate}
Then there can be a.s.\ at most one infinite component.
\end{thm}

In order to be able to use this result, we will first embed the metric generating tree into $\mathbb{Z}$ in a stationary (and ergodic) way. 
The embedding will be such that for each $r$, we have 
\begin{itemize}
\item[a.] any ball of radius $r$ will be represented by $N^r$ consecutive integers,
\item[b.]  the collection of balls of radius $r$ partitions $\mathbb{Z}$.
\end{itemize}
We first describe the construction rather loosely, and after that provide a formal construction. For ease of description, a collection of $m$ consecutive integers is called an {\em interval of length} $m$.

The ball of radius 1 containing 0, that is, $\ball_1(0)$ is chosen uniformly at random among all $N$ possible intervals of length $N$ containing the origin of $\mathbb{Z}$. Once we have
chosen this ball, all other balls of radius 1 are determined by requirements (a) and (b) above, although it is not yet clear at this point to exactly which
balls in $\Omega_N$ they correspond. To get an idea of this first step of the procedure, note that for $N=2$ there are only two possibilities, one of which is depicted here: 

\medskip
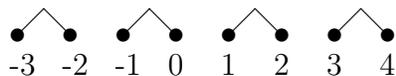
\begin{center}
\begin{figure}[!h]
\centering
\begin{picture}(150,30)
\put(3,15){\circle*{5}}
\put(23,15){\circle*{5}}
\put(43,15){\circle*{5}}
\put(63,15){\circle*{5}}
\put(83,15){\circle*{5}}
\put(103,15){\circle*{5}}
\put(123,15){\circle*{5}}
\put(143,15){\circle*{5}}
\put(3,15){\line(1,1){10}}
\put(23,15){\line(-1,1){10}}
\put(43,15){\line(1,1){10}}
\put(63,15){\line(-1,1){10}}
\put(83,15){\line(1,1){10}}
\put(103,15){\line(-1,1){10}}
\put(123,15){\line(1,1){10}}
\put(143,15){\line(-1,1){10}}
\put(0,0){-3}
\put(20,0){-2}
\put(40,0){-1}
\put(60,0){0}
\put(80,0){1}
\put(100,0){2}
\put(120,0){3}
\put(140,0){4}
\end{picture}
\label{fig:first}
\end{figure}
\end{center}
The other possibility is obtained by translating the edges over one unit to the right (or to the left, for that matter).

Next, we determine $\ball_2(0)$. The ball $\ball_2(0)$ is a union of $N$ balls of radius 1 and contains $\ball_1(0)$. There are $N$
possible ways to achieve this, keeping in mind that any ball of radius $2$ must - according to (a) above - be an interval of
length $N^2$. We now simply choose one of the $N$ possible ways to do this, with probability $1/N$ each. Once we have
chosen $\ball_2(0)$, all other balls of radius $2$ are determined for the same reason as before.
The following picture illustrates a possible choice for $\ball_2(0)$ given the choice of $\ball_1(0)$ made before.

\medskip
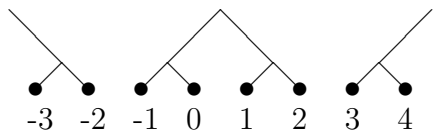
\begin{figure}[!h]
\centering
\begin{picture}(150,40)
\put(3,15){\circle*{5}}
\put(23,15){\circle*{5}}
\put(43,15){\circle*{5}}
\put(63,15){\circle*{5}}
\put(83,15){\circle*{5}}
\put(103,15){\circle*{5}}
\put(123,15){\circle*{5}}
\put(143,15){\circle*{5}}

\put(3,15){\line(1,1){10}}
\put(23,15){\line(-1,1){30}}
\put(43,15){\line(1,1){30}}
\put(63,15){\line(-1,1){10}}
\put(83,15){\line(1,1){10}}
\put(103,15){\line(-1,1){30}}
\put(123,15){\line(1,1){30}}
\put(143,15){\line(-1,1){10}}

\put(0,0){-3}
\put(20,0){-2}
\put(40,0){-1}
\put(60,0){0}
\put(80,0){1}
\put(100,0){2}
\put(120,0){3}
\put(140,0){4}

\end{picture}
\label{fig:second}
\end{figure}
We can continue this procedure as long as we wish, and in doing so we obtain a metric generating tree which is isomorphic
to the tree depicted in Figure 1. This last statement perhaps requires some reflection: one can see that this holds by first
identifying the two 0's in both graphs, and then build up the balls $\ball_r(0)$, $r=1,2,\ldots$, in that order. 

It is intuitively clear that this construction yields a stationary metric generating tree in the sense that the distribution of the stochastic
process which assigns to each pair $\{z,z'\}$ of points in $\mathbb{Z}$ the distance between them, is invariant under integer
translations. However, we would like to formalise the construction in such a way that not only stationarity follows
as an easy corollary, but we also obtain that the embedding of the metric generating tree is in fact ergodic with respect to translations. 

A possible formal construction is the following. Our probability space is the unit interval $[0,1]$ endowed with Lebesgue 
measure on its Borel sigma field. For $\gamma \in [0,1]$, let $\gamma=0.\gamma_1 \gamma_2\cdots$ be its $N$-adic expansion,
that is,
$$
\gamma=\sum_{n=1}^{\infty}\gamma_n N^{-n},
$$  
where we ignore those $\gamma$ for which the expansion is not unique - this is a set of Lebesgue measure zero anyway. In the
construction above, we saw that for each $r$, $\ball_{r-1}(0)$ can be seen as one of the balls of radius $r-1$ among the balls making
up $\ball_{r}(0)$. The metric generating tree corresponding to $\gamma \in [0,1]$ is obtained as follows. We let $\ball_r(0)$ be
such that $\ball_{r-1}(0)$ is the $(\gamma_r + 1)$-st ball in $\ball_{r}(0)$, counted from left to right. For instance, in the preceding
two figures with $N=2$, we have that $\gamma_1=1$ and $\gamma_2=0$. The map which assigns to each (apart from the exceptional
null set discussed before) $\gamma$ a metric generating tree is denoted by $\phi$. This map $\phi$ is invertible on a set
of full Lebesgue measure. 

It is clear that this construction formalises the
informal description given earlier. Furthermore, one can write down explicitly the transformation $S:[0,1] \to [0,1]$ 
which corresponds with the left-shift $T$ on the space of metric generating trees in the sense that $\phi \circ S = T \circ \phi$, hence $T=\phi S \phi^{-1}$. Indeed, a little reflection shows that $S$ can be described as follows: if $Y(\gamma)=\min\{n; \gamma_n \neq N-1\}$ then $S(\gamma)_k$ (that is, the $k$-th digit in $S(\gamma)$) is given by
$$
S(\gamma)_k =\left\{
\begin{array}{ll}
0 & \mbox{ if } k < Y(\gamma),\\
\gamma_k +1 & \mbox{ if } k=Y(\gamma),\\
\gamma_k & \mbox{ if } k > Y(\gamma).
\end{array}
\right.
$$
This transformation has been studied in the literature and goes by the name Kakutani - Von Neumann transformation, see
e.g.\ \cite{fried70}.
It is easy to check that Lebesgue measure is invariant under the action of $S$, and this immediately proves that the
construction of our random metric generating tree is stationary on $\mathbb{Z}$. 

\medskip\noindent
{\em Proof of Theorem \ref{theorem1a}.} 
The construction above shows that the metric generating tree can be embedded into $\mathbb{Z}$ in a stationary way. We claim that
this implies that the whole long range percolation process on the hierarchical lattice can be realised as a stationary
percolation process on $\mathbb{Z}$. To see this, we assign a uniformly-$[0,1]$ distributed random variable $U_e$ to
each edge $e$ in such a way that the collection is independent.

Given a realisation of the metric generating tree, we declare
edge $e$ to be open if $U_e \leq 1-\exp(-\alpha/\beta^{|e|})$, where $|e|$ denotes the length of $e$. This gives a realisation
of the percolation process with the correct distribution, and shows that we have embedded the full percolation process on
the hierarchical tree in a stationary way. Since every pair of vertices shares an edge, with positive probability, irrespective of the presence or absence of other edges, the positive finite energy condition is met and the result  follows.
\hfill{$\qed$}

\medskip\noindent
With a little more work one can also see that the
construction is in fact ergodic, that is, any event which is invariant under the shift on $\mathbb{Z}$ has probability 0 or 1.
To show this, we first show that Lebesgue measure on $[0,1]$ is ergodic with respect to $S$. 
This result is known, but
we give a simple (and new) proof for the convenience of the reader.

\begin{lem}
\label{eva}
Lebesgue measure on $[0,1]$ is ergodic with respect to $S$.
\end{lem}

\medskip\noindent
{\em Proof.}
Consider the first digit in each of $\gamma$, $S(\gamma)$, $S^2(\gamma), \ldots$. From the construction we immediately
conclude that the first digit follows the periodic pattern $0,1,\ldots, N-1,0,1,2,\ldots, N-1, \ldots$ starting at any number. Hence the first digit is just adding 1 modulo $N$. The second digit can only change when the first digit is an $N-1$, and then it
also changes according to adding 1 modulo $N$. In general the $k$-th digit can only change when the $(k-1)$-st digit is an $N-1$,
and then the change consists of adding 1 modulo $N$. It follows from these observations that the orbit of $\gamma$ under
the action of $S$ visits any $N$-adic interval $I_{m,k}= [kN^{-m}, (k+1)/N^{-m}]$ with frequency $N^{-m}$, that is,
\begin{equation}
\label{ergav}
\lim_{n \to \infty} \frac{1}{n} \sum_{j=0}^{n-1} {\bf 1}_{S(\gamma) \in I_{m,k}} = N^{-m},
\end{equation}
where ${\bf 1}_A$ denotes the indicator function of $A$.

Now consider the collection ${\cal M}$ of invariant probability measures for the transformation $S$. From the fact that 
Lebesgue measure preserves measure under $S$ we see that ${\cal M}$ is not empty. It is well known and easy to see
that the set ${\cal M}$ is convex, and that the ergodic measures are precisely the extremal points of ${\cal M}$. Since
${\cal M} \neq \emptyset$, this implies that there is at least one ergodic measure with respect to $S$.

Let $\nu$ be any ergodic measure with respect to $S$. It then follows from the ergodic theorem and (\ref{ergav}) that $\nu(I_{m,k}) = N^{-m}$ for any $m$ and $k=0, \ldots, 
N^m-1$. However, there is only one measure that satisfies this condition, namely Lebesgue measure on $[0,1]$. Hence $\nu$ must
be Lebesgue measure, which we already know is indeed invariant.
\hfill{$\qed$}

\begin{thm}\label{ergodic}
The embedding of our long range percolation process on the hierarchical lattice into $\mathbb{Z}$ is ergodic.
\end{thm}

\medskip\noindent
{\em Proof.}
From Lemma \ref{eva} it follows that the metric generating tree is embedded ergodically. Adding the i.i.d.\ random variables
$U_e$ as before does not destroy ergodicity, and the final configuration is a factor of this ergodic process and hence ergodic itself.
\hfill{$\qed$}

\section{Proof of Theorem \ref{theorem3}}
The proof consists of three steps:

\begin{enumerate}
\item For every constant $K>0$ the indicator function of the event that both $|\mathcal{C}(0)|=\infty$ and $|\mathcal{C}_n(0)| < K (\beta/N)^n$ converges a.s.\ to 0 as $n \to \infty$.
\item The fraction of the vertices in $\balln(0)$ which are in a cluster of size at least $K(\beta/N)^n$, converges a.s.\ to $\theta$ as $n \to \infty$.
\item Combine the previous two steps.
\end{enumerate}

\medskip\noindent
{\bf Step 1.}
We compute 
$$\mathbb{P}\left(|\mathcal{C}(0)|=\infty\Bigl||\{n \in \mathbb{N};|\mathcal{C}_n(0)| \leq K \left(\frac{\beta}{N}\right)^n\}| = \infty\right).$$
Let $n_1$ be the smallest $n$ for which $|\mathcal{C}_n(0)| \leq K (\beta/N)^n$, if $\mathcal{C}_{n_i}(0) \leftrightarrow  \overline{\mathcal{B}_{n_i}(0)}$, then $n_{i+1}$ is the smallest $n>n_i$ for which $\mathcal{C}_{n_i}(0) \not\leftrightarrow  \overline{\mathcal{B}_{n}(0)}$ and for which $|\mathcal{C}_n(0)| \leq K (\beta/N)^n$.\

Since $|\mathcal{C}_{n_i}(0)| \leq K (\beta/N)^{n_i}$, we have
\begin{equation}\label{inflimrep2}
\begin{array}{rcl}
\mathbb{P}\left(\mathcal{C}_{n_i}(0) \leftrightarrow \overline{\mathcal{B}_{n_i}(0)}\right) &
\leq & \mathbb{P} \left( \mathcal{C}_i(0) \leftrightarrow  \overline{\mathcal{B}_i(0)} \Bigl| |\mathcal{C}_i(0)| = \left\lfloor K \left(\frac{\beta}{N}\right)^{i} \right\rfloor \right)\\
\ & \leq & 1-\exp\left(-\alpha K \left(\frac{\beta}{N}\right)^{i}  \displaystyle\sum_{j=i+1}^{\infty} (N-1) \frac{N^{j-1}}{\beta^j}\right)\\
\ & = & 1-\exp\left(-\alpha K \frac{N-1}{\beta-N}\right).
\end{array}
\end{equation}
The right hand side is strictly less than 1 and is independent of $n_i$. So there will be an $n_i$ for which $\{\mathcal{C}_{n_i}(0) \not\leftrightarrow  \overline{\mathcal{B}_{n_i}(0)}\}$, and it follows that
$$\mathbb{P}\left(|\mathcal{C}(0)|=\infty \Bigl||\{n \in \mathbb{N};|\mathcal{C}_n(0)| \leq K \left(\frac{\beta}{N}\right)^n\}| = \infty\right) = 0.$$

\medskip\noindent
{\bf Step 2.} 
We use the random embedding of the hierarchical lattice in $\mathbb{Z}$, introduced in the previous section.
By Theorem \ref{ergodic} and the ergodic theorem we have, for every $k>0$, 
\begin{multline*}
\frac{1}{2 N^{n}+1} \sum_{x=-N^n}^{N^n} \mathbf{1}\left(\bigcap_{j=k}^{\infty}\left\lbrace  |\mathcal{C}_j(x)|>K\left(\frac{\beta}{N}\right)^j \right\rbrace  \right) \aston \\
\mathbb{P}\left(\bigcap_{j=k}^{\infty} \left\lbrace |\mathcal{C}_j(0)|>K\left(\frac{\beta}{N}\right)^j \right\rbrace \right).
\end{multline*}
From step 1 we know that this final probability increases to $\theta$ as $k \to \infty$, and it follows that
$$
A(n):=\frac{1}{2 N^n+1} \sum_{x=-N^n}^{N^n} \mathbf{1}\left(|\mathcal{C}_n(x)|>K\left(\frac{\beta}{N}\right)^n\right) \aston \theta.
$$
Note that the collection vertices $\{-N^n,-N^n+1,-N^n+2 \cdots, N^n\}$ contains the image under the embedding of the ball $\mathcal{B}_n(0)$ and this image contains a fraction $N^n/(2N^n+1)$ of those vertices.
Whether or not $|\mathcal{C}_n(x)|>K(\beta/N)^n$ is independent for vertices in different $n$-balls, so $$A_1(n):= \frac{1}{2 N^n+1} \sum_{x \in \mathcal{B}_n(0)} \mathbf{1}\left(|\mathcal{C}_n(x)|>K\left(\frac{\beta}{N}\right)^n\right)$$ and $A_2(n) := A(n) -A_1(n)$ are independent. Furthermore, $A_1(n)$ and $A_2(n)$ are bounded above by $1$ and have asymptotically the same mean. Combining these observations with 
$$
A_1(n) + A_2(n) = A(n) \aston \theta
$$ 
finishes the proof of Step 2.

\medskip\noindent
{\bf Step 3.} 
The strategy is to split those components in $\mathcal{B}_{n+1}(0)$ which are at least of size $K (\beta/N)^n$ into clusters roughly of size $K (\beta/N)^n$. Then we use those clusters as ``meta-vertices'' for an $N$-partite graph, in which meta-vertices in different $n$-balls are connected if the clusters they represent are connected by an edge of length $n+1$. Meta-vertices in the same $n$-ball never share an edge. We show that if we choose $K$ and $n$ large enough, then the largest component of the graph of meta-vertices contains a fraction of the meta-vertices close to 1, which shows that for large $n$, the fraction of vertices in the largest cluster of $\mathcal{B}_{n+1}(0)$ is close to $\theta$. We will be more precise now.

By step 2 we know that for every $K > 0$, every $\varepsilon>0$ and all large enough $n$, it holds that
$$\mathbb{P}\left(\left|\left\lbrace x \in \mathcal{B}_n(0);|\mathcal{C}_n(x)|>K\left(\frac{\beta}{N}\right)^n \right\rbrace \right|>(\theta-\varepsilon) N^n\right)>1-\varepsilon.$$
We now fix $\varepsilon$. The ball $\mathcal{B}_{n}(y)$ is said to be {\em good} if
$$\left|\left\lbrace x \in \mathcal{B}_n(0);|\mathcal{C}_n(x)|>K\left(\frac{\beta}{N}\right)^n \right\rbrace \right|>(\theta-\varepsilon) N^n.$$
We condition on the event that all $n$-balls in $\mathcal{B}_{n+1}(0)$ are good. The probability of this event is bounded below by $(1-\varepsilon)^N> 1-N \varepsilon$.

Now, for every good ball $\mathcal{B}_n(y)$, $y\in \Omega_N$, we make a partition of the set 
$$\left\lbrace x \in \mathcal{B}_n(y);|\mathcal{C}_n(x)|>K \left(\frac{\beta}{N}\right)^n \right\rbrace$$
 in ``meta-vertices''. For the moment we denote this set by $\mathcal{B}'_n(y)$. For $x \in \mathcal{B}'_n(y)$ we make a partition of $\mathcal{C}_n(x)$ in $\lfloor |\mathcal{C}_n(x)|/(\lceil K(\beta/N)^n \rceil)\rfloor$ sets, which all have size at least $\lceil K(\beta/N)^n \rceil$. Here $\lceil x \rceil := \inf\{n \in \mathbb{Z}; n \geq x\}$ is the ceiling of $x$ and $\lfloor x \rfloor := \sup\{n \in \mathbb{Z}; n \leq x\}$ is the floor of $x$.
The vertices that are not in such a cluster are ignored for the moment. Denote the collection of meta-vertices that contain vertices in $\mathcal{B}_{n+1}(0)$ by $\mathcal{V}_n$. We note that if $\mathcal{B}_{n}(y)$ is good and $K$ is large enough, then it contains at least $(\theta-\varepsilon) N^n/\lceil 2 K(\beta/N)^n\rceil \geq (\theta-\varepsilon) N^n/(3 K(\beta/N)^n)$ vertices.

We construct a new $N$-partite graph on $\mathcal{V}_n$ as follows. Let $\mathcal{V}_n$ be the vertex set and let $\mathcal{E}_n$ be the set of edges between those vertices.
This edge set is obtained as follows. Choose $\lceil K(\beta/N)^n \rceil$ original vertices from every meta-vertex $\mathcal{V}_n$. Chosing those vertices may be done in any way that is independent of the presence of edges of lenght $n+1$ or larger. Denote these sets by  $\mathcal{A}_n$.
The meta-vertices $x, y \in \mathcal{V}_n$ share an edge in $\mathcal{E}_n$, if there is at least 1 edge in the original graph that is shared by vertices that make up the sets in $\mathcal{A}_n$ corresponding to $x$ and $y$, and if the original vertices that make up $x$ and $y$ are at distance $n+1$ of each other. Otherwise there is no edge between the meta-vertices.

As observed before, the number of meta-vertices in $\mathcal{V}_n$ that consist of vertices from a good ball $\mathcal{B}_n(x)$, is at least $(\theta-\varepsilon) N^n/(3 K(\beta/N)^n)$. Since $\beta<N^2$, this quantity grows to $\infty$ as $n \to \infty$.
The expected degree of a vertex in $\mathcal{V}_n$ exceeds
$$\frac{(N-1) (\theta-\varepsilon) N^n}{3  K(\beta/N)^n}\left(1-\exp\left(-\alpha \beta^{-(n+1)} (K)^2  \left(\frac{\beta}{N}\right)^{2n}\right)\right),$$
which is larger than $\lambda:=(\tilde{N}-1) (\theta-\varepsilon) \alpha K/(6 \beta)$, for all large enough $n$. This holds for every $K >0$, and therefore the expected degree can be chosen to be arbitrary large.

This $N$-partite graph falls within the class of inhomogeneous random graphs of Bollob{\'a}s, Janson and Riordan \cite{bollobas05}.
The degree of every meta-vertex is asymptotically Poisson distributed, with mean bounded below by $\lambda$ and we know that the (unique) largest component of such an $N$-partite graph contains with high probability (in the limit for $n \to \infty$) a fraction $\rho$ of the meta-vertices, where $\rho$ is the largest solution of $$1-\rho = e^{-\lambda \rho}.$$ By tuning $K$, $\lambda$ can be chosen arbitrary large and $\rho$ can be taken such that $\rho >1 -\varepsilon$.
So, for every $\varepsilon>0$ and large enough $n$ the graph $(\mathcal{V}_n,\mathcal{E}_n)$ contains a unique giant component, which contains a fraction $(1 - \varepsilon) N$ of the vertices in $\mathcal{V}_n$, with probability at least $1-\varepsilon$.

Since we have conditioned on the event that all $n$-balls in $\mathcal{B}_{n+1}(0)$ are good, the fraction of vertices in $\mathcal{B}_{n+1}(0)$, that are part of vertices in $\mathcal{V}_n$ is bounded below by $\theta - 2\varepsilon$. (The factor 2, is due to the fact that the sizes of different meta-vertices differ at most a factor 2). Therefore, with the same conditioning, the largest cluster in $\mathcal{B}_{n+1}(0)$ is at least of size $$(\rho - \varepsilon)(\theta-2\varepsilon) N^n > (1 - 2 \varepsilon)(\theta-2\varepsilon) N^n$$ with probability exceeding
$1-\varepsilon$. Now multiplying by the probability that all $n$-balls in $\mathcal{B}_{n+1}(0)$ are good, gives that the probability that the largest cluster in $\mathcal{B}_{n+1}(0)$ is at least of size $(1 - 2 \varepsilon)(\theta-2\varepsilon) N^n$ is bounded below by $(1-\varepsilon)(1-N \varepsilon)$.
By chosing $\varepsilon' \in (0,\varepsilon/\max(4,N+1))$, we obtain that
$\mathbb{P}(|\bar{\mathcal{C}}_n(0)|>(\theta-\varepsilon') N^n)$ is at least $1-\varepsilon'$ and this finishes the proof. 
\hfill{$\qed$}

\medskip\noindent
\textbf{Remark} We realize that it is possible to prove the statement of Step 2 by using the strong law of large numbers. If we do this, then it is only a small step from the proof of Theorem \ref{theorem3} to a proof of Theorem \ref{theorem1a}. However, we think that the proof presented in the previous section contains some valuable ideas and therefore should be included in this paper.

\section{Proof of Theorem \ref{theorem2}}

Continuity proofs of percolation functions typically split into separate proofs for left and right continuity, one of which typically follows
from standard arguments. In this case, continuity from the right in $\alpha$ and continuity from the left in $\beta$ are the easy parts:

\begin{lem}\label{leftcontlem}
$\theta(\alpha,\beta)$ is continuous from the right in $\alpha >0$ and continuous from the left in $\beta>0$.
\end{lem}

\medskip\noindent
{\em Proof.}
For $\alpha>0$ and $\beta \leq N$, $\theta(\alpha,\beta)=1$, so the statement of the lemma holds in that domain.
Note that $$ \{\mathcal{C}_i(0) \leftrightarrow \overline{\mathcal{B}_i(0)}\} \subset
\{\mathcal{C}_{i-1}(0) \leftrightarrow \overline{\mathcal{B}_{i-1}(0)}\}$$ and
\begin{equation}\label{inflimrep}
\{|\mathcal{C}(0)|=\infty\} = \cap_{i=0}^{\infty} \{\mathcal{C}_i(0) \leftrightarrow \overline{\mathcal{B}_i(0)}\}.
\end{equation}
A straightforward computation yields that for $\beta > N$,
\begin{equation}\label{inflimrep3}
\mathbb{P}\left(\mathcal{C}_i(0) \leftrightarrow \overline{\mathcal{B}_i(0)}\right) = \mathbb{E}\left(1-\exp\left(-\alpha |\mathcal{C}_i(0)|  \frac{N-1}{\beta-N}\left(\frac{N}{\beta}\right)^i\right)\right).
\end{equation}
Since $|\mathcal{C}_i(0)|$ depends on the state of only finitely many edges, this expectation is continuous in $\alpha$ and $\beta$. In particular, $\mathbb{P}(\mathcal{C}_i(0) \leftrightarrow \overline{\mathcal{B}_i(0)})$ is continuous from the left in $\beta$ for $\beta>N$ and therefore it is continuous from the left for $\beta>0$ (in fact, it is also continuous from the right).
Furthermore, $\mathbb{P}(\mathcal{C}_i(0) \leftrightarrow \overline{\mathcal{B}_i(0)})$ is increasing in $\alpha$, decreasing in $\beta$ and decreasing in $i$.
Because a decreasing limit of increasing (resp.\ decreasing) functions, which are continuous from the right (resp.\ left)
is continuous from the right (resp. left), the statement follows.
\hfill{$\qed$}

\medskip
In order to prove that $\theta(\alpha,\beta)$ is continuous from the left in $\alpha >0$ and continuous from the right in $\beta>0$, we use a renormalisation argument. Fix $\alpha> 0$ and $N \leq \beta < N^2$. To get insight in the argument we first (falsely) assume that for given $\varepsilon>0$ and large enough finite $k$, there is a $\delta>0$ such that,
$$\mathbb{P}_{\alpha-\delta,\beta+\delta}(|\mathcal{C}^m_k(0)| > (\theta(\alpha,\beta) - \varepsilon) N^k) = 1.$$
(Although the assumption is false, we can get this probability arbitrary close to 1, by choosing $k$ large enough, $\delta$ small enough  and using Theorem \ref{theorem3}.)

Now we use renormalisation. The balls of radius $k$ are considered as vertices of $\Omega_N$ which we call ``meta-vertices". If two vertices in the original model have distance $k+l$, then the meta-vertices in which they are contained are at distance $l$. Vertices in the new model are connected if and only if the largest clusters in the 
original $k$-balls, represented by these vertices, are connected by an edge. The new model is again a percolation model on $\Omega_N$.

Let $x$ and $y$ be meta-vertices, at distance $l$ of each other. Define, for $\delta>0$ small,  
$$
\alpha' := (\alpha-\delta)((\theta(\alpha,\beta) - \varepsilon))^2 \frac{N^{2k}}{(\beta+\delta)^{k}}.
$$
Given the states of all other edges, the (conditional) probability that $x$ and $y$ are connected to each other is always bounded below by
$$
1-\exp(-(\alpha-\delta)((\theta(\alpha,\beta) - \varepsilon) N^k)^2 (\beta+\delta)^{-(k+l)})
$$
and by the choice of $\alpha'$, this is just
$$
1-\exp(-\alpha' (\beta+\delta)^{-l}).
$$
Hence, the renormalized model stochastically dominates the percolation model with parameters $\alpha'$ and $\beta+\delta$.

Since $N^2/(\beta+\delta)>1$, $\alpha'$ can be chosen arbitrary large by choosing $k$ large. In particular it can be chosen such that $\theta(\alpha',\beta+\delta) > 1-\varepsilon$, (by the second remark after the proof of Theorem \ref{theorem1}). It follows that for large enough $k$
\begin{eqnarray*}
\mathbb{P}{}_{\alpha-\delta,\beta+\delta}\left(|\mathcal{C}| = \infty\right) &\geq &
\theta(\alpha', \beta + \delta) \mathbb{P}{}_{\alpha-\delta,\beta+\delta}(0 \in \mathcal{C}^m_k(0))\\
&\geq &(1-\varepsilon)(\theta(\alpha,\beta) - \varepsilon) \\
&\geq & \theta(\alpha,\beta) - 2 \varepsilon.
\end{eqnarray*}

The only problem is that we have incorrectly assumed that
$$
\mathbb{P}{}_{\alpha-\delta,\beta+\delta}\left(|\mathcal{C}^m_k(0)| > (\theta(\alpha,\beta) - \varepsilon) N^k\right) = 1,
$$
and we will deal with this problem now. We need the notion of mixed percolation (cf.~\cite{chayesschonmann}). Mixed percolation involves independently removing vertices, together with all of its adjacent edges.
Formally, the measure $\mathbb{P}^{mixed}_{\alpha,\beta,\gamma}$ is constructed as follows. A vertex in $\Omega_N$ is open with probability $1-\gamma$, independently of the states (open or closed) of the other vertices in $\Omega_N$. If $x$ and $y$ are both open vertices then they share an edge with probability $1-\exp(-\alpha/\beta^{d(x,y)})$, independently of which other edges are present.

\begin{lem}\label{mixedlem} Let $\beta > N$.
For every $\varepsilon>0$, there exists $\gamma>0$, such that  $$\mathbb{P}{}_{\alpha,\beta}(|\mathcal{C}(0)|=\infty) \leq \mathbb{P}{}^{mixed}_{\alpha(1+\varepsilon),\beta,\gamma}(|\mathcal{C}(0)| =\infty).$$
\end{lem}

Before giving the proof of this result, we show how it can be used to prove Theorem \ref{theorem2}. The following lemma suffices.

\begin{lem}\label{contlemma3}
If $\theta(\alpha,\beta)>0$, then  for every $\epsilon \in (0,\theta(\alpha,\beta))$, there exists $\delta>0$ such that $\theta(\alpha-\delta,\beta+\delta)> \theta(\alpha,\beta)-\epsilon$.
\end{lem}

\medskip\noindent
{\em Proof.}
Fix $\alpha$, $\beta$ and $\epsilon>0$.
Let $\alpha'$ be such that
$$
\theta(\alpha',(N^2+\beta)/2) > 1- \varepsilon/3,
$$
which is possible by Theorem \ref{theorem1}(b) and the second remark after its proof. Furthermore, let $\gamma \in (0,\epsilon/3)$ be such that
$$\theta(\alpha',(N^2+\beta)/2) \leq  \mathbb{P}{}^{mixed}_{2 \alpha',(N^2+\beta)/2,\gamma}(|\mathcal{C}(0)| =\infty),$$
which is possible by Lemma \ref{mixedlem}.
Let $K$ be such that the following conditions are satisfied:
\begin{enumerate}
\item $\alpha(\theta(\alpha,\beta)-\epsilon/2)2 (N^{2}/\beta)^K > 3 \alpha'$,
\item $\mathbb{P}_{\alpha,\beta}(|\mathcal{C}_K(0)|>(\theta(\alpha,\beta)-\epsilon/3)N^K)> 1-\gamma/2$,
\end{enumerate}
which are possible by respectively $N^2>\beta$ and Theorem \ref{theorem3}.
Finally, let $\delta>0$ be such that  $\delta < \min(\alpha/3,(N^2-\beta)/2)$ and
\begin{equation}\label{mixedincont}
\mathbb{P}{}_{\alpha-\delta,\beta+\delta}(|\mathcal{C}^m_K(0)|>(\theta(\alpha,\beta)-\epsilon/3)N^K)> 1-\gamma,
\end{equation}
which is possible by the continuity of the probability in $\alpha$ and $\beta$ for finite $K$.

We say that the ball $\mathcal{B}_K(x)$ is {\em good} if $\mathcal{C}^m_K(x)$ has size at least $(\theta(\alpha,\beta)-\epsilon/3) N^K$. Delete all vertices that are in a ball of diameter $K$ which is not good and also all vertices that are not in the largest cluster of good balls.
As above, we interpret the remaining components as the vertices of  the hierarchical lattice of order $N$ in which vertices are independently deleted with probability at most $\gamma$, by (\ref{mixedincont}).
Remaining clusters in the original graph, of which the vertices are at distance $K+l$, are connected by at least one edge with probability at least
$$
1-\exp(-(\alpha -\delta)(\theta(\alpha,\beta)-\hat{\epsilon}/3)2 N^{2K}\beta^{-(K+l)})> 1-\exp(-2\alpha'\beta^{-l}),
$$
irrespective of the existence or absence of other connections. Here we have used that $\alpha -\delta> 2 \alpha/3$. Hence the rescaled process stochastically dominates a mixed percolation process with parameters
$2 \alpha'$, $\beta$ and $\gamma$.

Now note that by exchangability
\begin{eqnarray*}
\ & \ & \mathbb{P}{}_{\alpha-\delta,\beta+\delta}(|\mathcal{C}_K(0)| \geq  (\theta(\alpha,\beta)-\hat{\epsilon}/3)N^K) \\
&\geq &(\theta(\alpha,\beta)-\hat{\epsilon}/3) \mathbb{P}_{\alpha-\delta,\beta+\delta}(|\mathcal{C}^m_K(0)|>(\theta(\alpha,\beta)-\hat{\epsilon}/3)N^K) \\
&\geq & (\theta(\alpha,\beta)-\hat{\epsilon}/3)(1-\gamma).
\end{eqnarray*}
Furthermore, conditioned on $0$ being in the largest cluster of a good ball, the probability that 0 is in an infinite cluster if the parameters are $\alpha-\delta$ and $\beta+\delta$ is larger than $1- \epsilon/3$. Combining these observations and $\gamma < \epsilon/3$ gives
$$\theta(\alpha-\delta,\beta+\delta) > (\theta(\alpha,\beta)-\hat{\epsilon}/3)(1-\epsilon/3)2 > \theta(\alpha,\beta)-\epsilon.$$
\hfill{$\qed$}

\medskip\noindent
It remains to prove Lemma \ref{mixedlem}.
Before giving a proof of this lemma we define a directed long-range percolation model and relate this to the undirected model. In the directed version, vertices in $\Omega_N$ are open with probability $1-\gamma$. If vertex $x$ is open, then a directed edge from $x$ to $y$ is present with probability $1-\exp\left(-\alpha \beta^{-d(x,y)}\right)$. Conditioned on the states of the vertices (open or closed) the presence or absence of an edge is independent of the presence or absence of other edges. The corresponding measure we denote by $\hat{\mathbb{P}}^{mixed}_{\alpha,\beta,\gamma}$. The set of vertices which can be reached by a path from vertex $x$ is denoted by $\hat{\mathcal{C}}(x)$. Note that in th directed model, the presence of a path from $x$ to $y$ does not necessarily imply that there exists a path from $y$ to $x$. We define the directed version of the original (not mixed) measure, $\hat{\mathbb{P}}_{\alpha,\beta}$, in a similar way and note that $\hat{\mathbb{P}}_{\alpha,\beta} = \hat{\mathbb{P}}^{mixed}_{\alpha,\beta,0}$.
Standard arguments (see e.g.\ \cite{cox88,meester09}) can be used to show that 
$$\mathbb{P}^{mixed}_{\alpha,\beta, \gamma}(|\mathcal{C}(0)|=\infty) = \hat{\mathbb{P}}^{mixed}_{\alpha,\beta, \gamma}(|\hat{\mathcal{C}}(0)|=\infty).$$

\medskip\noindent
{\em Proof of Lemma \ref{mixedlem}.}
The directed mixed percolation graph with parameters $\alpha$, $\beta$ and $\gamma$ can be obtained as follows (the ordinary model can be obtained by taking $\gamma=0$). We assign i.i.d.\ random variables $X_x$ to the vertices $x \in \Omega_N$, all Poisson distributed with parameter $\alpha (N-1)/(\beta-N)$. We construct a directed multi-graph (a graph in which multiple edges between two vertices in the same direction are allowed). Vertices are open with probability $1-\gamma$, independently of each other.  If $x$ is open, then $X_x$ directed edges start at $x$. The endpoints of these edges are independently chosen from $\Omega_N \setminus x$, and a vertex at distance $r$ of $x$ is chosen with probability $(\beta-N)(N-1)^{-1}\beta^{-r}$. If $x$ is closed, then no edges start at $x$.
We obtain the original directed graph by replacing the collection of all edges from $x$ to $y$ (if there is at least one) by a single edge from $x$ to $y$, for all $x,y \in \Omega_N$.

Let $Z_1$ be a Poisson distributed random variable with parameter $\alpha (N-1)/(\beta-N)$. Furthermore, let $Z_2 = Y_1Y_2$, where $Y_1$ is equal to 1 with probability $1-\gamma$ and equal to 0 with probability $\gamma$, and where $Y_2$ is independent of $Y_1$ and Poisson distributed with parameter $\alpha (1+\varepsilon) (N-1)/(\beta-N)$.
For the ordinary percolation model the number of edges starting at $x$ in the multigraph is distributed as $Z_1$, while for the mixed percolation model, the number of edges starting at $x$ is distributed as $Z_2$.
It is now easy to check that for $\varepsilon >0$ there is a $\gamma> 0$, such that $\mathbb{P}(Z_1 = 0) = \mathbb{P}(Z_2 = 0)$ and for this $\gamma$ and all $k>0$ we have,
$$\mathbb{P}(Z_2>k|Z_2 > 0) = \mathbb{P}(Y_2>k|Y_2 > 0) > \mathbb{P}(Z_1>k|Z_1 > 0).$$ The statement of Lemma \ref{mixedlem} now follows by a straightforward coupling argument.
\hfill{$\qed$}

\section{Proof of Theorem \ref{theorem4}}

The proof of Theorem \ref{theorem4} is split into separate proofs of continuity from the right and from the left of  $\alpha_{c}(\beta)$.

\begin{lem}
$\alpha_{c}(\beta)$ is strictly increasing on $\beta \in (N,N^2)$ and  continuous from the right on $\beta \in (0,N^2)$.
\end{lem}

\medskip\noindent
{\em Proof}
Theorem \ref{theorem2} implies that $\theta(\alpha_{c}(\beta),\beta)=0$ for $\beta > N$. By observing that for every $\varepsilon>0$, $\mathbb{P}_{(1+\varepsilon)\alpha,(1+\varepsilon)\beta}$ is stochastically dominated by $\mathbb{P}_{\alpha,\beta}$, we deduce that $\alpha_{c}(\beta(1+\varepsilon)) \geq (1+\varepsilon) \alpha_{c}(\beta)$. Since by Theorem \ref{theorem1}, $\alpha_{c}(\beta) > 0$ for $\beta \in (N,N^2)$, this gives that $\alpha_{c}(\beta)$ is strictly increasing on $(N,N^2)$.

In order to prove continuity from the right, we note that for all $\delta>0$, $\theta(\alpha_{c}(\beta)+\delta,\beta)>0$ by definition. By the continuity of $\theta(\alpha,\beta)$ we obtain that there exist $\varepsilon>0$, such that $\theta(\alpha_{c}(\beta)+\delta,\beta+\varepsilon)>0$ and therefore $\alpha_{c}(\beta+\varepsilon)< \alpha_{c}(\beta)+\delta$. This together with $\alpha_{c}(\beta+\varepsilon) > \alpha_{c}(\beta)$  completes the proof.
\hfill{$\qed$}

\begin{lem}\label{leftcontcr}
$\alpha_{c}(\beta)$ is  continuous from the left for $\beta \in (0,N^2)$.
\end{lem}

To prove this, we use \cite{aizenman87}. In that paper it is shown that for long-range percolation on $\mathbb{Z}^d$,
\begin{equation}\label{equalcrits}
\inf\{\alpha:\theta(\alpha,\beta) > 0\} = \sup\{\alpha:\mathbb{E}_{\alpha,\beta}(|\mathcal{C}(0)|) < \infty\}.
\end{equation}
Inspection of the proof of this result yields that this proof also works on the hierarchical lattice.
Now use the following lemma.
\begin{lem}\label{contexp}
Let $\alpha>0$ and $\beta>N$. If $\mathbb{E}_{\alpha,\beta}(|\mathcal{C}(0)|) < \infty$, then there exist $\varepsilon>0$ such that
$\mathbb{E}_{\alpha,\beta(1-\varepsilon)}(|\mathcal{C}(0)|) < \infty$.
\end{lem}

\medskip\noindent
{\em Proof of Lemma \ref{leftcontcr} (given Lemma \ref{contexp}).}
For every $\varepsilon>0$, we have that $$\mathbb{P}_{\alpha_{c}(\beta),\beta(1-\varepsilon)}(|\mathcal{C}(0)| = \infty)> 0,$$ by the strict increase of $\alpha_{c}(\beta)$. This implies that $\mathbb{E}_{\alpha_{c}(\beta),\beta(1-\varepsilon)}(|\mathcal{C}(0)|) = \infty$, for every $\varepsilon>0$ and therefore $\mathbb{E}_{\alpha_{c}(\beta),\beta}(|\mathcal{C}(0)|) = \infty$.

Furthermore, by (\ref{equalcrits}) we know that for every $\delta>0$, $\mathbb{E}_{\alpha_{c}(\beta)-\delta,\beta}(|\mathcal{C}(0)|) < \infty$. Therefore, there exist $\varepsilon>0$ such that $\mathbb{E}_{(\alpha_{c}(\beta)-\delta),\beta(1-\varepsilon)}(|\mathcal{C}(0)|) < \infty$, which implies that for all $\delta>0$, there exist $\varepsilon>0$ such that $\alpha_{c}(\beta(1-\varepsilon))> \alpha_{c}(\beta)-\delta$. This  together with  $\alpha_{c}(\beta(1-\varepsilon)) < \alpha_{c}(\beta)$ gives continuity from the left of $\alpha_{c}(\beta)$.
\hfill{$\qed$}

\medskip\noindent
{\em Proof of Lemma \ref{contexp}.}
Assign independent uniform$(0,1)$ random variables to all pairs of vertices in $\Omega_N$. The random variable assigned to the pair $(x,y)$ is denoted by $U(x,y)$. We say that $x$ and $y$ share an edge for the parameters $\alpha$ and $\beta$ if $U(x,y)<1-\exp(-\alpha \beta^{-d(x,y)})$. This construction provides a coupling for long-range percolation models with different values of $\alpha$ and $\beta$. Define $\mathcal{C}(x;\alpha,\beta)$ as the cluster of vertices that can be reached by paths if the parameters are $\alpha$ and $\beta$. 

Assume that $a := \mathbb{E}_{\alpha,\beta}(|\mathcal{C}(0)|) < \infty$ and take $\varepsilon >0$ small enough (we will
see later exactly how small). Define $\mathcal{A}_0(0) := \mathcal{C}(0;\alpha,\beta)$. In an inductive fashion,
let $\mathcal{A}_{i+1}'(0)$ be the set of vertices not in $\cup_{j=0}^{i} \mathcal{A}_{j}(0)$ that can be reached from $\mathcal{A}_i(0)$ by crossing an edge present for the parameters $\alpha$ and $\beta(1-\varepsilon)$.
Note that $\mathcal{A}'_i(0) \subset \mathcal{A}_i(0)$ and that 
by construction, $\mathcal{C}(0;\alpha,\beta(1-\varepsilon)) = \cup_{i=1}^{\infty} \mathcal{A}_i(0)$. The next step in the proof is to bound  $\mathbb{E}(|\mathcal{A}_i(0)|)$. By definition $\mathbb{E}(|\mathcal{A}_0(0)|)=a$. Since the graph is transitive, for every $x \in \mathcal{A}_i(0)$ the expected size of the set
$$\{y \in \Omega_N \setminus \cup_{j=0}^{i}\mathcal{A}_j(0); U(x,y) < 1-\exp(-\alpha ((1-\varepsilon)\beta)^{-d(x,y)})\}$$
is bounded above by
\begin{eqnarray*}
b & :=& \sum_{y \in \Omega_N} \mathbb{P} \Bigl( U(0,y) < 1-\exp(-\alpha ((1-\varepsilon)\beta )^{-d(0,y)}) \Bigl|
U(0,y) >1-\exp[-\alpha \beta^{-d(0,y)}] \Bigr)\\
&=& \sum_{y \in \Omega_N} \frac{\mathbb{P}(1-\exp(-\alpha \beta^{-d(0,y)})   < U(0,y) < 1-\exp(-\alpha ((1-\varepsilon)\beta)^{-d(0,y)})}{\mathbb{P}(1-\exp(-\alpha \beta^{-d(0,y)})   < U(0,y))}\\
& =&  \sum_{i = 1}^{\infty}  (N-1)N^{i-1} \frac{\exp(-\alpha \beta^{-i}) -\exp (-\alpha \beta^{-i}(1-\varepsilon)^{-i})}{\exp(-\alpha \beta^{-i})} \\
&=& \sum_{i = 1}^{\infty}  (N-1)N^{i-1} (1 -\exp (-\alpha \beta^{-i}((1-\varepsilon)^{-i} -1))) \\
&\leq & \sum_{i = 1}^{\infty} (N-1)N^{i-1} \alpha \beta^{-i}[(1-\varepsilon)^{-i} -1] \\
&=& (N-1)\alpha \frac{(1-\varepsilon)}{(1-\varepsilon) \beta -N} - (N-1)\alpha \frac{1}{\beta -N} \\
&=& \frac{\alpha \varepsilon (N-1)N}{(\beta(1-\varepsilon)-N)(\beta-N)},
\end{eqnarray*}
which converges to 0, if $\varepsilon \searrow 0$. Therefore, we can choose $\varepsilon>0$ so small that $b < a^{-1}$.
Note that $\mathbb{E}(|\mathcal{A}'_{i+1}(0)|) \leq b \mathbb{E}(|\mathcal{A}_{i}(0)|)$, and because of the transitivity of the graph,
$\mathbb{E}(|\mathcal{A}_{i+1}(0)|) \leq a \mathbb{E}(|\mathcal{A}'_{i+1}(0)|)$.
So, we have $$\mathbb{E}(|\mathcal{A}_{i}(0)|) \leq (ab)^i \mathbb{E}(|\mathcal{A}_{0}(0)|).$$

Since $\mathcal{C}(0;\alpha,\beta(1-\varepsilon)) = \cup_{i=1}^{\infty} \mathcal{A}_i(0)$ and $ab<1$, we have
$$\mathcal{C}(0;\alpha,\beta(1-\varepsilon)) \leq \sum_{i=0}^{\infty} a (ab)^i = \frac{a}{1-ab} < \infty.$$
\hfill{$\qed$}

\medskip\noindent
{\em Proof of Theorem \ref{theorem4}.}
The only things left to prove is that $\alpha_{c}(\beta) \nearrow \infty$ for $\beta \nearrow N^2$, but this follows immediately from Lemma \ref{contexp}, equality (\ref{equalcrits}) and observing that $\alpha_{c}(N^2) = \infty$.
\hfill{$\qed$}

\section{Possible generalizations}

Possible generalizations of the model considered in this paper include:

\begin{enumerate}
\item {\em Randomness in the hierarchical lattice.} The hierarchical lattice $\Omega_N$ is generated by a $N$-regular tree. An interesting question is how randomness in the underlying tree (or induced random metric) affects the percolation process on the resulting lattice.
One possibility is that the metric generating tree, is a Galton-Watson tree. 
Analysis of long-range percolation on such random hierarchical structures is not a trivial extension of the analysis in this paper, since in renormalisation schemes, one has to take care of all kinds of dependencies of the sizes of balls of given diameter.

\item {\em More general connection function $p(k)$.} In this paper we focused on $p_k=1-\exp\left(-\frac{\alpha}{\beta^k}\right)$. What are necessary and sufficient conditions on $g(k)$ so that when $p_k=1-\exp\left(-\alpha g(k)\right)$ we have $0 < \alpha_{c} < \infty$?

\item {\em Random cluster models.} We only consider independent percolation on the hierarchical lattice. We did not try to incorporate Random cluster (or Fortuin-Kasteleyn) model \cite{Grim06} yet. Some work has already been done for the Ising model on the hierarchical lattice \cite{hara01}.
\end{enumerate}

\medskip\noindent
{\bf Acknowledgments} We thank Odo Diekmann for suggesting to study this percolation model. Research of VK was supported by NWO grant NWO-CLS Project 635.100.002. Research of PT was supported by a grant from Riksbankens jubileumsfond.

\bibliographystyle{plain}
\bibliography{hierpercarx}

\begin{thebibliography}{10}

\bibitem{aizenman87}
M.~Aizenman and D.J. Barsky.
\newblock {Sharpness of the phase transition in percolation models}.
\newblock {\em Comm. Math. Phys.}, 108(3):489--526, 1987.

\bibitem{aizenman86}
M.~Aizenman and C.M. Newman.
\newblock Discontinuity of the percolation density in one-dimensional {$1/\vert
  x- y\vert \sp 2$} percolation models.
\newblock {\em Comm. Math. Phys.}, 107(4):611--647, 1986.

\bibitem{AtSw09}
S.M. Athreya and J.M. Swart.
\newblock Survival of contact processes on the hierarchical group.
\newblock {\em Probab. Theory Relat. Fields}, 2009.
\newblock In Press: DOI 10.1007/s00440-009-0214-x.

\bibitem{berger02}
N.~Berger.
\newblock Transience, recurrence and critical behavior for long-range
  percolation.
\newblock {\em Comm. Math. Phys.}, 226(3):531--558, 2002.

\bibitem{biskup04}
M.~Biskup.
\newblock On the scaling of the chemical distance in long-range percolation
  models.
\newblock {\em Ann. Probab.}, 32(4):2938--2977, 2004.

\bibitem{biskup09}
M.~Biskup.
\newblock Graph diameter in long-range percolation.
\newblock preprint at arXiv:math/0406379v2 [math.PR], 2009.

\bibitem{bollobas05}
B.~Bollob{\'a}s, S.~Janson, and O.~Riordan.
\newblock {The phase transition in inhomogeneous random graphs}.
\newblock {\em Random Structures and Algorithms}, 31(1-2):3--122, 2007.

\bibitem{chayesschonmann}
L.~Chayes and R.H. Schonmann.
\newblock {Mixed percolation as a bridge between site and bond percolation}.
\newblock {\em Ann.~Appl.~Probab.}, 10(4):1182--1196, 2000.

\bibitem{cox88}
J.T. Cox and R.~Durrett.
\newblock {Limit theorems for the spread of epidemics and forest fires}.
\newblock {\em Stochastic processes and their applications}, 30(2):171--191,
  1988.

\bibitem{dawson07}
D.A. Dawson and L.G. Gorostiza.
\newblock Percolation in a hierarchical random graph.
\newblock {\em Commun. Stoch. Anal.}, 1(1):29--47, 2007.

\bibitem{deme90}
F.M. Dekking and R.W.J. Meester.
\newblock On the structure of mandelbrot's percolation process and other random
  cantor sets.
\newblock {\em J. Stat. Phys.}, 58:1109--1126, 1990.

\bibitem{fried70}
N.A. Friedman.
\newblock {\em {Introduction to ergodic theory}}.
\newblock Van Nostrand Reinhold, New York, 1970.

\bibitem{gandolfi92}
A.~Gandolfi, M.S. Keane, and C.M. Newman.
\newblock Uniqueness of the infinite component in a random graph with
  applications to percolation and spin glasses.
\newblock {\em Probab. Theory Related Fields}, 92(4):511--527, 1992.

\bibitem{Grim06}
G.R. Grimmett.
\newblock {\em The random cluster model}, volume 333 of {\em Grundlehren der
  Mathematischen Wissenschaften [Fundamental Principles of Mathematical
  Sciences]}.
\newblock Springer-Verlag, Berlin, 2006.

\bibitem{hara01}
T.~Hara, T.~Hattori, and H.~Watanabe.
\newblock {Triviality of hierarchical Ising model in four dimensions}.
\newblock {\em Commun. Math. Phys.}, 220(1):13--40, 2001.

\bibitem{meester09}
R.~Meester and P.~Trapman.
\newblock {Bounding basic characteristics of spatial epidemics with a new
  percolation model}.
\newblock preprint at arXiv:0812.4353v2 [math.PR].

\bibitem{newmanChuck86}
C.M. Newman and L.S. Schulman.
\newblock {One dimensional $1/\vert j-i\vert \sp s$ percolation models: The
  existence of a transition for {$s\leq 2$}}.
\newblock {\em Communications in Mathematical Physics}, 104(4):547--571, 1986.

\bibitem{schulman83}
L.S. Schulman.
\newblock {Long range percolation in one dimension}.
\newblock {\em Journal of Physics A: Mathematical and General}, 16:L639--L641,
  1983.

\bibitem{trapman10}
P.~Trapman.
\newblock {The growth of the infinite long-range percolation cluster}.
\newblock to appear in Annals of Probability, preprint at arXiv:0901.0661v2
  [math.PR], 2010.

\bibitem{vanrooij78}
A.C.M. Van~Rooij.
\newblock {Non-Archimedean functional analysis, vol. 51 of Monographs and
  Textbooks in Pure and Applied Math}, 1978.

\end{thebibliography}

\end{document}